\documentclass[10 pt]{amsart}
\usepackage[applemac]{inputenc}
\usepackage{amsaddr}

\usepackage[colorlinks=true,linkcolor=red,citecolor=blue,urlcolor=blue,hypertexnames=false]{hyperref}
\usepackage{bookmark}
\usepackage{amsthm,thmtools,amssymb,amsmath,amscd}
\usepackage{xcolor}

\usepackage{fancyhdr}
\usepackage{esint}
\usepackage{enumerate}

\usepackage{pictexwd,dcpic}

\swapnumbers
\declaretheorem[name=Theorem,numberwithin=section]{thm}
\declaretheorem[name=Remark,style=remark,sibling=thm]{rem}
\declaretheorem[name=Lemma,sibling=thm]{lemma}

\declaretheorem[name=Example,style=remark,sibling=thm]{example}

\numberwithin{equation}{section}

\usepackage{cleveref}
\crefname{lemma}{Lemma}{Lemmata}
\crefname{prop}{Proposition}{Propositions}
\crefname{thm}{Theorem}{Theorems}
\crefname{cor}{Corollary}{Corollaries}
\crefname{defn}{Definition}{Definitions}
\crefname{example}{Example}{Examples}
\crefname{rem}{Remark}{Remarks}
\crefname{assum}{Assumption}{Assumptions}
\crefname{nota}{Notation}{Notation}

\usepackage{autonum}

\newcommand{\ti}{\tilde}

\newcommand{\cn}{\colon}
\newcommand{\sub}{\subset}

\newcommand{\bbR}{\mathbb{R}}

\newcommand{\bbS}{\mathbb{S}}

\newcommand{\8}{\infty}

\newcommand{\al}{\alpha}
\newcommand{\be}{\beta}
\newcommand{\g}{\gamma}

\newcommand{\ka}{\kappa}
\newcommand{\la}{\lambda}

\newcommand{\s}{\sigma}

\newcommand{\Om}{\Omega}
\newcommand{\D}{\Delta}
\newcommand{\G}{\Gamma}

\newcommand{\cL}{\mathcal{L}}

\newcommand{\cW}{\mathcal{W}}

\newcommand{\cD}{\mathcal{D}}

\newcommand{\cV}{\mathcal{V}}

\newcommand{\cK}{\mathcal{K}}
\newcommand{\cR}{\mathcal{R}}
\newcommand{\cE}{\mathcal{E}}

\newcommand{\del}{\partial}
\newcommand{\n}{\nabla}

\newcommand{\fa}{\forall}

\newcommand{\ip}[2]{\left\langle #1,#2 \right\rangle}
\newcommand{\fr}[2]{\frac{#1}{#2}}
\newcommand{\x}{\times}

\DeclareMathOperator{\id}{id}
\DeclareMathOperator{\pr}{pr}

\DeclareMathOperator{\graph}{graph}

\DeclareMathOperator{\tr}{tr}
\DeclareMathOperator{\Rm}{Rm}

\DeclareMathOperator{\EV}{EV}

\newcommand{\pf}[1]{\begin{proof}{\parskip\baselineskip{ #1}} \end{proof}}
\newcommand{\eq}[1]{\begin{equation}\begin{alignedat}{2} #1 \end{alignedat}\end{equation}}

\newcommand{\br}[1]{\left(#1\right)}
\newcommand{\abs}[1]{\lvert #1\rvert}
\newcommand{\enum}[1]{\begin{enumerate}[(i)] #1 \end{enumerate}}
\newcommand{\enu}[1]{\begin{enumerate}[(a)] #1 \end{enumerate}}

\newcommand{\Ra}{\Rightarrow}
\newcommand{\ra}{\rightarrow}

\newcommand{\mt}{\mapsto}

\newcommand{\mrm}{\mathrm}

\newcommand{\hp}{\hphantom}
\newcommand{\q}{\quad}

\begin{document}
\title[Concavity of solutions to degenerate elliptic equations]{Concavity of solutions to degenerate elliptic equations on the sphere}
\date{\today}
\thanks{The work of JS is funded by the "Deutsche Forschungsgemeinschaft" (DFG, German research foundation); Project "Quermassintegral preserving local curvature flows"; Grant number SCHE 1879/3-1.}

\author[M. Langford]{Mat Langford}
\address{Department of Mathematics, University of Tennessee Knoxville, TN 37996, USA}
\email{\href{mailto:mlangford@utk.edu}{mlangford@utk.edu}}

\author[J. Scheuer]{Julian Scheuer}
\address{Department of Mathematics, Columbia University
New York, NY 10027, USA}
\email{\href{mailto:julian.scheuer@math.uni-freiburg.de}{julian.scheuer@math.uni-freiburg.de}}

\begin{abstract}
We prove the concavity of classical solutions to a wide class of degenerate elliptic differential equations on strictly convex domains of the unit sphere. The proof employs a suitable two-point maximum principle, a technique which originates in works of Korevaar, Kawohl and Kennington for equations on Euclidean domains. We emphasize that no differentiability of the differential operator is needed, but only some monotonicity and concavity properties. 
\end{abstract}

\maketitle

\section{Introduction}
This paper is about classical solutions to fully nonlinear degenerate equations of the form
\eq{\label{Equ}f(\abs{\n u},-\n^{2} u)=b(\cdot,u,\abs{\n u})}
on a domain $\Om\sub \bbS^{n}$, $n\geq 2$, of the unit sphere whose closure $\bar\Om$ is geodesically convex. (We recall that a subset $C$ of $\bbS^n$ is geodesically convex if no two points of $C$ are antipodal and the unique minimizing geodesic segment between any two points of $C$ lies entirely in $C$.) Here
\eq{\n u=(Du)^{\sharp}\;\;\text{and}\;\; \n^{2}u=(D^{2}u)^{\sharp},}
where $D$ is the Levi-Civita connection of the round metric $\s$ on $\bbS^{n}$ and the $\sharp$-operator is the usual index raising operator with respect to $\s$, so that $\n u$ is the gradient vector field and $\n^{2}u$ the Hessian endomorphism field on $\Om$. 

We want to prove concavity of solutions $u\in C^{2}(\Om)\cap C^{1}(\bar\Om)$ to \eqref{Equ} satisfying suitable boundary conditions.
The question of concavity of solutions to elliptic equations on Euclidean domains  is widely studied, in particular in connection with Laplacian eigenvalues. For example, Brascamp and Lieb \cite{BrascampLieb:08/1976} proved that the first Dirichlet eigenfunction of the Laplacian over a convex domain is log-concave. A recent highlight in this area was, with many partial results preceding it, the proof of the fundamental gap conjecture by Andrews and Clutterbuck \cite{AndrewsClutterbuck:07/2011}.
Using a refined log-concavity estimate for the first eigenfunction, they obtained a sharp estimate for the difference between the first and second Dirichlet eigenvalues (the ``fundamental gap''). An analogous result on the sphere was obtained by Seto, Wang and Wei \cite{SetoWangWei:/2019}.

Our main tool will be the ``two-point maximum principle'' approach introduced by Korevaar \cite{Korevaar:/1983}. The idea is to estimate the function
\eq{\label{TPF-1}\ti Z(\la,x,y)=u(\la x+(1-\la)y)-(\la u(x)+(1-\la)u(y),\q 0\leq \la\leq 1,}
using a first and second derivative test from below by zero, where $u$ is the solution to a certain quasilinear equation with some structure properties. Refinements of this approach for $C^{2}$-solutions are due to Kennington \cite{Kennington:/1985} and Kawohl \cite{Kawohl:/1985}. Using regularization, Sakaguchi \cite{Sakaguchi:/1987} has transferred this method to solutions of certain $p$-Laplace equations, the solutions of which may be non-smooth. 
Similar extensions to viscosity solutions are the content of \cite{AlvarezLasryLions:03/1997} and \cite{Juutinen:10/2010}. The argument was used to obtain proofs of differential Harnack inequalities for hypersurface flows in \cite{BourniLangford:12/2019}.

To our knowledge, there is no such concavity maximum principle for non-Euclid\-ean domains. The reason for this seems to arise from a major technical obstruction: in order to make the application of the maximum principle to \eqref{TPF-1} work, it is necessary to differentiate a family of geodesics with respect to their end-points. This leads to the consideration of a family of Jacobi fields and their variations. In order to exploit the first order conditions, one needs detailed knowledge about the nonlinear Jacobi fields along the $\ti Z$-minimizing geodesic. Even worse, the second order condition contains a further space derivative of the Jacobi equation and in general it seems impossible to extract the correct sign, which is required to apply the maximum principle, from these conditions.

The purpose of this paper is to initiate the study of concavity maximum principles for degenerate elliptic equations on a Riemannian manifold and we start this journey on the unit sphere. Here we have managed to deal with the pretty well known Jacobi fields and their maybe not so well known derivatives.
One first key observation is that, for technical reasons, the full \emph{three}-parameter function $\ti Z$ is no longer suitable. The concavity of a solution is equivalent to the concavity of
\eq{t\mt u(\g(t))}
for every geodesic 
\eq{\g\cn [-1,1]\ra \bar\Om,}
 cf. \cite{Udriste:/1994}. By a well known result due to Jensen \cite{Jensen:/1906}, for continuous functions this is equivalent to asking this function to be {\it{midpoint-concave}}, i.e.
\eq{u(\g(0))\geq \fr 12 \br{u(\g(-1))+u(\g(1))}.}

Hence it will suffice to prove that the function
\eq{\label{Z}Z\cn \bar\Om\x\bar\Om&\ra \bbR\\
			(x,y)&\mt u\br{\g_{x,y}(0)}-\fr 12(u(\g_{x,y}(-1))+u(\g_{x,y}(1)))
}
is non-negative, where $\g_{x,y}$ is the unique minimizing geodesic from $x$ to $y$. The mid-point approach has also been taken for equations on Euclidean domains in \cite{Juutinen:10/2010,Kawohl:/1985}. 

Before we state our main theorem, we need to introduce some more objects. We denote by $(\cV,\cL)$ the category of $n$-dimensional real vector spaces with linear transformations as its objects. A function $f$ on that category is called {\it{isotropic}} if it acts on the subclass $\cE$ of endomorphisms in a $\mrm{GL}_{n}$-invariant way:
\eq{\label{ISO1}\fa \cW\in \cE~\fa V\in \cL_{\mrm{inv}}\cn f(\cW)=f(V^{-1}\circ \cW\circ V),}
where $\cL_{\mrm{inv}}$ denotes the class of all invertible linear maps, understood to map to the correct space such that \eqref{ISO1} makes sense.
From this property it is clear that the action of $f$ on real diagonalizable endomorphism is only determined through the action on the ordered set of eigenvalues $(\ka_{1},\dots,\ka_{n})$, 
\eq{\label{ISO2}f(\cW)=:\ti f(\EV(\cW))=\ti f(\ka_{1},\dots,\ka_{n}).}
The properties we will require from $f$ will in general only be satisfied if $f$ is restricted to a subclass of $\cE$ and the best way to describe the domain of definition of $f$ is via $\ti f$. Therefore we suppose the domain of definition for $\ti f$, $\G\sub\bbR^{n}$, is a symmetric, open and convex cone which contains the positive cone
\eq{\G_{+}=\{\ka\in \bbR^{n}\cn \ka_{i}>0~\fa 1\leq i\leq n\}.}
Denote by $\mrm{Sym}^{n\x n}(\bbR)$ the vector space of symmetric $(n\x n)$-matrices. Then there exists an open domain of definition $\cD_{\G}\sub \mrm{Sym}^{n\x n}(\bbR)$ for $f$ consisting of matrices with eigenvalues in $\G$, such that \eqref{ISO2} holds for all $\cW\in \cD_{\G}$. 

Abusing notation, the orthogonal invariance of $f$ allows us to understand $\cD_{\G}$ to be a subset of any other space of $\bar g$-self-adjoint endomorphisms of an arbitrary $n$-dimensional inner product space $(E,\bar g)$. 
Here is our main theorem:
 
\begin{thm}\label{thm:main}
Let $n\geq 2$, $\Om\sub\bbS^{n}$ and suppose $\bar\Om$ to be geodesically convex.  Let $\G\sub\bbR^{n}$ be a symmetric, open and convex cone containing $\G_{+}$ and let $u\in C^{2}(\Om)\cap C^{1}(\bar\Om)$ be a solution to
\eq{f(\abs{\n u},-\n^{2}u)=b(\cdot,u,\abs{\n u})\q \text{in}~\Om.}
Suppose the function $f\cn [0,\8)\x  \cD_{\bar\G}\ra \bbR$ has the following properties:
\enum{
\item For every $p\in [0,\8)$, $f(p,\cdot)$ is an isotropic function on ${\cD}_{\bar\G},$
\item $f$ is increasing in the first variable, 
\eq{p\leq q\q\Ra\q f(p,\cdot)\leq f(q,\cdot), }
\item increasing in the second variable, 
\eq{\ka_{i}\leq \la_{i}\q\fa 1\leq i\leq n\q\Ra\q f(\cdot,\mrm{diag}(\ka_{1},\dots,\ka_{n}))\leq f(\cdot,\mrm{diag}(\la_{1},\dots,\la_{n})),}
\item convex in the second variable, 
\eq{f(\cdot,\tfrac 12(A+B))\leq \tfrac 12 (f(\cdot,A)+f(\cdot,B))\q\fa A,B\in \cD_{\bar\G}.}
}
Suppose the function $b\cn \Om\x \bbR\x [0,\8)\ra \bbR$ is
\enu{
\item decreasing in the third variable,
\item jointly concave in the first two variables:
	\eq{b\br{\g_{x,y}(0),\tfrac 12(u(x)+u(y)),p}\geq \tfrac 12 b(x,u(x),p)+\tfrac 12 b(y,u(y),p),}
\item strictly decreasing in the second variable.
}
If furthermore for all $(x,y)\in \del\Om\x\bar\Om$ there holds
\eq{Du_{x}(\dot{\g}_{x,y}(-1))-Du_{y}(\dot{\g}_{x,y}(1))>0,}
  then $u$ is concave.
\end{thm}

\begin{rem}
 The boundary condition is slightly stronger than the ones in \cite{Kennington:/1985,Korevaar:/1983}. This stems from our technical restriction that due to the nonlinear ambient geometry we have to use the notion of mid-point concavity. Hence we can not vary a boundary point without varying another point too. Geometrically it says that at every point $x\in \del\Om$, the totally geodesic hypersurface tangent to $\graph(u)$ at $(x,u(x))$ must lie above $\graph(u)$ {\it{and}} at every point $y\in \Om$ the totally geodesic hypersurface tangent to $\graph(u)$ at $(y,u(y))$ must lie above $\graph(u_{|\del\Om})$, and one of these relations must be strict.  
\end{rem}

\begin{example}
Let us give a large class of examples of operators to which this theorem applies. Given a convex and increasing function $\psi\cn \bbR\ra \bbR$, define
\eq{\ti f(\ka_{1},\dots,\ka_{n})=\sum_{i=1}^{n}\psi(\ka_{i}).}
Then the operator function
\eq{f(\cW)=\tr(\psi(\cW))}
associated to $\ti f$ is convex and increasing. Here $\psi$ is also understood to be the operator function. Explicit examples are 
\eq{f(-\n^{2}u)=-\D u, \q f(-\n^{2}u)=\tr(\exp(-\n^{2}u)),}
where the latter is the matrix exponential.  
\end{example}

For the proof of \Cref{thm:main} we need an easy result from linear algebra and some tedious calculations of Jacobi fields and their derivatives. We address these issues in the next two sections, before we can proceed to the proof of the theorem.

\section{A result from linear algebra}

The nonlinear ambient space will introduce an algebraic problem.
In order deal with it, we need the following lemma from linear algebra. It should be standard, but we did not find a proper reference, so let us give the proof for the reader's convenience.

\begin{lemma}\label{Courant}
Let $n\geq 1$ and $\cW$ and $V$ be symmetric $(n\x n)$-matrices, such that $\cW$ is non-negative definite and $V$ defines a bijective contraction map, i.e.
\eq{\abs{V(x)}\leq \abs{x}.}
Then the spectra of the matrices $\cW$ and $V\circ \cW\circ V$ are ordered in the sense
\eq{\la_{i}\leq \ka_{i},}
where $\la_{1}\leq \dots\leq \la_{n}$ are the eigenvalues of $V\circ\cW\circ V$ and $\ka_{1}\leq \dots\leq\ka_{n}$ are those of $\cW$.
The reverse inequality holds if $V$ is an expansion.
\end{lemma}

\pf{Let $\cR_{A}$ denote the Rayleigh quotient of a matrix,
\eq{\cR_{A}(x)=\fr{\ip{Ax}{x}}{\abs{x}^{2}}.}
Then
\eq{\cR_{V\circ\cW\circ V}(x)=\fr{\ip{(\cW\circ V)x}{Vx}}{\abs{x}^{2}}\leq \fr{\ip{(\cW\circ V)x}{Vx}}{\abs{Vx}^{2}}=\cR_{\cW}(Vx),}
where we used the assumptions on $\cW$ and $V$.
From the Courant min-max principle we get
\eq{\la_{i}&=\min_{U}\left\{\max_{x}\{\cR_{V\circ \cW\circ V}(x)\cn x\in U\}\cn \mrm{dim}(U)=i\right\}\\
		&\leq \min_{U}\left\{\max_{x}\{\cR_{\cW}(Vx)\cn x\in U\}\cn \mrm{dim}(U)=i\right\}\\
		&= \min_{U}\left\{\max_{x}\{\cR_{\cW}(Vx)\cn x\in V^{-1}(U)\}\cn \mrm{dim}(U)=i\right\}\\
		&= \min_{U}\left\{\max_{x}\{\cR_{\cW}(y)\cn y\in U\}\cn \mrm{dim}(U)=i\right\}\\
		&=\ka_{i}.}
To prove the expansion case, just reverse all inequalities.
}

\section{Jacobi fields}
We have to differentiate \eqref{Z} twice and in this section we collect the geometric ingredients of this process. The derivatives of $\g_{x,y}$ will be Jacobi fields and variations of them and we will make extensive use of the Jacobi equation.
Define
\eq{\G\cn \Om\x\Om\x[-1,1]&\ra \Om\\
					(x,y,t)&\mt \g_{x,y}(t).}
For the moment choose arbitrary coordinates $(x^{i},y^{i})$ in $\Om\x\Om$.
The vector fields
\eq{J_{x^{i}}(x,y,\cdot)=\fr{\del \G}{\del x^{i}}(x,y,\cdot),\q J_{y^{i}}(x,y,\cdot)=\fr{\del \G}{\del y^{i}}(x,y,\cdot)}
are, for given $x,y$, Jacobi fields along $\g_{x,y}$ with	
\eq{J_{x^{i}}(x,y,-1)=\del_{x^{i}},\q J_{x^{i}}(x,y,1)=0,\q J_{y^{i}}(x,y,-1)=0,\q J_{y^{i}}(x,y,1)=\del_{y^{i}},}
where we identify $\del_{x^{i}}$ with its pushforward under the identity map $\G(\cdot,y,-1)$ and likewise for $\del_{y^{i}}$.				
We will have to deal with the second derivatives of $\G$ which are given by (neglecting the distinction between $D_{\del_{x^{j}}}$ and $D_{\fr{\del\G}{\del x^{j}}}$, since later we will work in Riemannian normal coordinates)
\eq{K_{x^{i}x^{j}}=D_{\del_{x^{j}}}\fr{\del\G}{\del x^{i}},\q K_{y^{i}y^{j}}=D_{\del_{y^{j}}}\fr{\del\G}{\del y^{i}},\q K_{x^{i}y^{j}}=D_{\del_{y^{j}}}\fr{\del\G}{\del x^{i}},\q K_{y^{i}x^{j}}=D_{\del_{x^{j}}}\fr{\del\G}{\del y^{i}},}
Those satisfy a differentiated Jacobi equation, which	we derive in the following.
We use the following convention for the Riemannian curvature tensor.
\eq{\Rm(X,Y)Z=D_{X}D_{Y}Z-D_{Y}D_{X}Z-D_{[X,Y]}Z,}
i.e. on $\Om\sub\bbS^{n}$ we have
\eq{\label{Riem}\Rm(X,Y)Z=\s(Y,Z)X-\s(X,Z)Y,}
where we recall that $\s$ denotes the round metric.
With this convention the Jacobi equation, which for example $J_{x^{i}}$ satisfies for fixed $x,y$, is
\eq{\label{Jac}0=D_{t}^{2}J_{x^{i}}+\Rm(J_{x^{i}},\dot{\G})\dot{\G}=D_{t}^{2}J_{x^{i}}+\abs{\dot\G}^{2}J_{x^{i}}-\s(\dot\G,J_{x^{i}})\dot\G,}
see e.g. \cite[Ch.~10]{Lee:/1997}. The same holds for $J_{y^{i}}$ and now we have to differentiate this equation with respect to $x^{j}$ and $y^{j}$.
We consider it to be easier to follow those calculations if we do not explicitly put in the specific form of the Riemann tensor just yet. However we will already use $D\Rm=0$.
 For $K_{x^{i}x^{j}}$ we calculate:

\begin{lemma}
There holds
\eq{\label{Ev-K}
0&=D_{t}^{2}K_{x^{i}x^{j}}+\Rm(K_{x^{i}x^{j}},\dot\G)\dot\G+2\Rm(J_{x^{i}},\dot\G)D_{t}J_{x^{j}}+2\Rm(J_{x^{j}},\dot\G)D_{t}J_{x^{i}}}
and similarly for $K_{y^{i}y^{j}}$, $K_{x^{i}y^{j}}$ and $K_{y^{i}x^{j}}$.
\end{lemma}
\pf{
We only prove this for $K_{x^{i}x^{j}}$ since the other identities are obtained in a similar manner. Differentiating the Jacobi equation and applying the first Bianchi identity (in the final step) yields
\eq{0&=D_{\del_{x^{j}}}\br{D_{t}^{2}J_{x^{i}}+\Rm(J_{x^{i}},\dot{\G})\dot{\G}}\\
	&=D_{t}D_{\del_{x^{j}}}D_{t}J_{x^{i}}+\Rm(J_{x^{j}},\dot \G)D_{t}J_{x^{i}}+\Rm(K_{x^{i}x^{j}},\dot\G)\dot\G+\Rm(J_{x^{i}},D_{t}J_{x^{j}})\dot\G\\
	&\hp{=}+\Rm(J_{x^{i}},\dot\G)D_{t}J_{x^{j}}\\
	&=D_{t}\br{D_{t}D_{\del_{x^{j}}}J_{x^{i}}+\Rm(J_{x^{j}},\dot\G)J_{x^{i}}}+\Rm(J_{x^{j}},\dot \G)D_{t}J_{x^{i}}+\Rm(K_{x^{i}x^{j}},\dot\G)\dot\G\\
	&\hp{=}+\Rm(J_{x^{i}},D_{t}J_{x^{j}})\dot\G+\Rm(J_{x^{i}},\dot\G)D_{t}J_{x^{j}}\\
	&=D_{t}^{2}K_{x^{i}x^{j}}+\Rm(D_{t}J_{x^{j}},\dot\G)J_{x^{i}}+2\Rm(J_{x^{j}},\dot\G)D_{t}J_{x^{i}}+\Rm(K_{x^{i}x^{j}},\dot\G)\dot\G\\
	&\hp{=}+\Rm(J_{x^{i}},D_{t}J_{x^{j}})\dot\G+\Rm(J_{x^{i}},\dot\G)D_{t}J_{x^{j}}\\
	&=D_{t}^{2}K_{x^{i}x^{j}}+\Rm(K_{x^{i}x^{j}},\dot\G)\dot\G+2\Rm(J_{x^{i}},\dot\G)D_{t}J_{x^{j}}+2\Rm(J_{x^{j}},\dot\G)D_{t}J_{x^{i}}
}
as claimed.
}				

In order to get more information about these quantities, first we note that at $t=\pm1$ all $K$ are zero. This is due to the fact that
\eq{\G(\cdot,\cdot,-1)=\pr_{x},\q \G(\cdot,\cdot,1)=\pr_{y}}
and the second covariant derivative of an isometry is always zero.
Also it is in order to pick the right coordinates. In the following, latin indices range from $0$ to $n-1$ and greek indices range from $1$ to $n-1$.
 At a given point $x_{0}$ pick Riemannian normal coordinates, where with a rotation we arrange
\eq{\del_{x^{0}}(x_{0})=\fr{\dot{\g}}{\abs{\dot\g}}(x_{0}),}
while at $y_{0}$ we do the same. Since we have Riemannian normal coordinates at $x_{0}$ and $y_{0}$, the $2n$ vectors $(\del_{x^{i}},\del_{y^{i}})$ form an orthonormal basis at $(x_{0},y_{0})$, while $\del_{y^{0}}$ is the parallel transport of $\del_{x^{0}}$. Furthermore, the Christoffel symbols vanish at $(x_{0},y_{0})$ in the product space. After a possible orthogonal transformation at $y_{0}$ and reordering of basis vectors, there exists a set of vector fields $(E_{i}(t))_{t\in[-1,1]}$ parallel along $\g_{x_{0},y_{0}}$, such that
\eq{E_{i}(-1)=\del_{x^{i}},\q E_{i}(1)=\del_{y^{i}}.}

\begin{lemma}\label{K}
In the above constructed coordinates there hold:
\enum{		
\item $K_{x^{0}x^{0}}=K_{y^{0}y^{0}}=K_{x^{0}y^{0}}=0.$
\item $K_{x^{0}x^{\al}}+K_{y^{0}y^{\al}}+K_{x^{0}y^{\al}}+K_{y^{0}x^{\al}}=0.$
\item The quantities
\eq{\cK_{2}=\br{K_{x^{\al}x^{\be}}+K_{y^{\al}y^{\be}}+ K_{x^{\al}y^{\be}}+ K_{y^{\al}x^{\be}}}\xi^{\al}\xi^{\be}}
and
\eq{\cK_{3}=\br{K_{x^{\al}x^{\be}}+K_{y^{\al}y^{\be}}- K_{x^{\al}y^{\be}}- K_{y^{\al}x^{\be}}}\xi^{\al}\xi^{\be}}
 where $\xi$ is any parallel vector field, satisfy
\eq{\label{K-2}0&=D_{t}^{2}\cK_{2}+\Rm(\cK_{2},\dot\G)\dot\G+4\Rm(J_{x^{\al}}+ J_{y^{\al}},\dot\G)(D_{t}J_{x^{\be}}+D_{t}J_{y^{\be}})\xi^{\al}\xi^{\be},}
respectively
\eq{\label{K-3}0&=D_{t}^{2}\cK_{3}+\Rm(\cK_{3},\dot\G)\dot\G+4\Rm(J_{x^{\al}}- J_{y^{\al}},\dot\G)(D_{t}J_{x^{\be}}- D_{t}J_{y^{\be}})\xi^{\al}\xi^{\be}.}
}
\end{lemma}
\pf{
(i)~Since $x_{0}$ and $y_{0}$ are not conjugate points, the Jacobi fields are uniquely determined by their boundary values. Hence, suppressing $(x_{0},y_{0})$, 
\eq{J_{x^{0}}(t)=\fr 12 (1-t)E_{0}(t),\q J_{y^{0}}(t)=\fr 12 (1+t) E_{0}(t).}
However, $E_{0}$ itself is the tangent to the normalized geodesic. Plugging this into \eqref{Ev-K} yields that those $K$ are Jacobi fields themselves, though with vanishing boundary values. Hence they are identically zero.

(ii)~For any $1\leq \al\leq n-1$ define
\eq{\cK_{1}=K_{x^{0}x^{\al}}+K_{y^{0}y^{\al}}+K_{x^{0}y^{\al}}+K_{y^{0}x^{\al}}.}
Then $\cK_{1}$ solves
\eq{0&=D_{t}^{2}\cK_{1}+\Rm(\cK_{1},\dot\G)\dot\G+2\br{\Rm(J_{x^{\al}},\dot\G)D_{t}J_{x^{0}}+\Rm(J_{y^{\al}},\dot\G)D_{t}J_{y^{0}}\right.\\
	&\hp{=}\left.+\Rm(J_{y^{\al}},\dot\G)D_{t}J_{x^{0}}+\Rm(J_{x^{\al}},\dot\G)D_{t}J_{y^{0}}}\\
	&=D_{t}^{2}\cK_{1}+\Rm(\cK_{1},\dot\G)\dot\G,}
where we have used
\eq{J_{x^{0}}+J_{y^{0}}=E_{0},\q D_{t}J_{x^{0}}=-\fr 12 E_{0}=-D_{t}J_{y^{0}}.}
Hence $\cK_{1}$ is, as a Jacobi field with vanishing boundary values, identically zero.

(iii)~ We add up all four equations of the form \eqref{Ev-K} with the appropriate sign. We treat $\cK_{2}$ and $\cK_{3}$ simultaneously by using $\pm$ in the appropriate places. We get
\eq{
0&=D_{t}^{2}K_{x^{\al}x^{\be}}+\Rm(K_{x^{\al}x^{\be}},\dot\G)\dot\G+2\Rm(J_{x^{\al}},\dot\G)D_{t}J_{x^{\be}}+2\Rm(J_{x^{\be}},\dot\G)D_{t}J_{x^{\al}}\\
	&\hp{=}+D_{t}^{2}K_{y^{\al}y^{\be}}+\Rm(K_{y^{\al}y^{\be}},\dot\G)\dot\G+2\Rm(J_{y^{\al}},\dot\G)D_{t}J_{y^{\be}}+2\Rm(J_{y^{\be}},\dot\G)D_{t}J_{y^{\al}}\\
	&\hp{=}\pm D_{t}^{2}K_{x^{\al}y^{\be}}\pm \Rm(K_{x^{\al}y^{\be}},\dot\G)\dot\G\pm 2\Rm(J_{x^{\al}},\dot\G)D_{t}J_{y^{\be}}\pm 2\Rm(J_{y^{\be}},\dot\G)D_{t}J_{x^{\al}}\\
	&\hp{=}\pm D_{t}^{2}K_{y^{\al}x^{\be}}\pm \Rm(K_{y^{\al}x^{\be}},\dot\G)\dot\G\pm 2\Rm(J_{y^{\al}},\dot\G)D_{t}J_{x^{\be}}\pm 2\Rm(J_{x^{\be}},\dot\G)D_{t}J_{y^{\al}}.
	}
Applying this to $\xi^{\al}\xi^{\be}$ gives
\eq{0&=D_{t}^{2}\cK_{2/3}+\Rm(\cK_{2/3},\dot\G)\dot\G+4\Rm(J_{x^{\al}},\dot\G)D_{t}J_{x^{\be}}\xi^{\al}\xi^{\be}\\
	&\hp{=}+4\Rm(J_{y^{\al}},\dot\G)D_{t}J_{y^{\be}}\xi^{\al}\xi^{\be}\pm 4\Rm(J_{x^{\al}},\dot\G)D_{t}J_{y^{\be}}\xi^{\al}\xi^{\be}\\
	&\hp{=}\pm 4\Rm(J_{y^{\al}},\dot\G)D_{t}J_{x^{\be}}\xi^{\al}\xi^{\be}\\
	&=D_{t}^{2}\cK_{2/3}+\Rm(\cK_{2/3},\dot\G)\dot\G+4\Rm(J_{x^{\al}}\pm J_{y^{\al}},\dot\G)D_{t}J_{x^{\be}}\xi^{\al}\xi^{\be}\\
	&\hp{=}+4\Rm(J_{y^{\al}}\pm J_{x^{\al}},\dot\G)D_{t}J_{y^{\be}}\xi^{\al}\xi^{\be}.}
A quick check of both cases reveals that $\cK_{2}$ and $\cK_{3}$ satisfy the desired equations.
}

\begin{lemma}\label{lem:Q=0}
The quantities $\cK_{2}$ and $\cK_{3}$ from \cref{K} both vanish at $t=0$.
\end{lemma}

\pf{
From \eqref{Jac} we obtain
\eq{J_{x^{\al}}(t)=v_{\al}(1-t)E_{\al}(t),\q J_{y^{\al}}(t)=v_{\al}(1+t)E_{\al}(t)}
with a function $v_{\al}$ that satisfies 
\eq{\label{v}\ddot{v}_{\al}=-\abs{\dot\G}^{2}v_{\al},\q v_{\al}(0)=0,\q v_{\al}(2)=1.}
 Hence, from multiplying \eqref{K-2} and \eqref{K-3} with some $E_{\be}$, we see that $\ip{\cK_{2/3}}{E_{\be}}E_{\be}$ is a Jacobi field with vanishing boundary values and hence zero.
Thus $\cK_{2}$ and $\cK_{3}$ are both multiples of $\dot\G$ and again from \eqref{K-2} and \eqref{K-3} we get
\eq{0&=\del_{t}^{2}\ip{\cK_{2/3}}{\dot\G}-4\Rm(J_{x^{\al}}\pm J_{y^{\al}},\dot\G,\dot\G,D_{t}J_{x^{\be}}\pm D_{t}J_{y^{\be}})\xi^{\al}\xi^{\be}\\
	&=-4(v_{\al}(1-t)\pm v_{\al}(1+t))\del_{t}(v_{\be}(1-t)\pm v_{\be}(1+t))\Rm(E_{\al},\dot\G,\dot\G,E_{\be})\xi^{\al}\xi^{\be}\\
	&\hp{=}+\del_{t}^{2}\ip{\cK_{2}}{\dot\G}.}
In this sum, the contribution of terms involving $\al\neq \be$ is zero, because in this case the curvature term is zero. So we get
\eq{0&=-2\abs{\dot{\G}}^{2}\del_{t}(v_{\al}(1-t)\pm v_{\al}(1+t))^{2}\xi^{\al}\xi^{\al}+\del_{t}^{2}\ip{\cK_{2}}{\dot\G}\\
	&=: -\dot g(t)+\ddot{h}(t).}  
We solve this equation by integration with respect to the boundary values $h(-1)=h(1)=0$ and evaluate $h$ at $t=0$.
In general we get with some constants $a$ and $b$,
\eq{\dot h(t)=\int_{-1}^{t}\dot g+a=g(t)-g(-1)+a}
and
\eq{h(t)=\int_{-1}^{t}g-(1+t)g(-1)+(1+t)a.}
Since $h(1)=0$, we find
\eq{a=g(-1)-\fr 12 \int_{-1}^{1}g}
and hence we have
\eq{h(t)=\int_{-1}^{t}g-\fr{1+t}{2} \int_{-1}^{1}g.}
Due to the symmetry of the function $g$ around $0$ we get $h(0)=0$ and the proof is complete.
}

\section{Proof of \texorpdfstring{\cref{thm:main}}{Theorem 1.1}}

To prove \Cref{thm:main}, we need to prove that the two-point function
\eq{Z(x,y)=u(\g_{x,y}(0))-\fr 12 (u(x)+u(y))}
defined in $\bar\Om\x\bar\Om$ is non-negative at minimal points. Therefore we may assume that $x\neq y$. If one of the minimizing points, say $x$, is at the boundary, then we prove that by moving $x$ and $y$ towards each other at the same speed (i.e. while fixing the mid-point $z$), we could decrease the value of $Z$. To see this, consider the function 
\eq{c(t)=u(z)-\fr 12(u(\g_{x,y}(t-1))+u(\g_{x,y}(1-t))).}
There holds
\eq{-2\dot{c}(0)=Du_{x}(\dot{\g}_{x,y}(-1))+Du_{y}(\dot{\g}_{x,y}(1))>0.}
Hence $\dot{c}(0)$ is strictly negative and $Z$ must achieve its minimum at interior points.

Let $z=\g_{x,y}(0)$.
The first order conditions on $Z$ are (suppressing the spatial arguments for $J$)
\eq{\label{1st order}0= Du_{z}(J_{x^{i}}(0))-\fr 12 \del_{i}u_x=v_i(1)\del_{{i}}u_{z}-\fr 12 \del_{i}u_x}
and
\eq{0=D u_{z}(J_{y^{i}}(0))-\fr 12 \del_{i}u_y=v_i(1)\del_{{i}}u_{z}-\fr 12 \del_{i}u_y,}
with 
\eq{v_{0}(1)=\fr 12,\q v_{\al}(1)=:v(1),}
and the latter solves \eqref{v} and hence
\eq{v(t)=\fr{\sin(\abs{\dot{\G}}t)}{\sin(2\abs{\dot\G})}.}
We use the bases
\eq{E_{i}(-1)=\del_{x^{i}},\q E_{i}(0),\q E_{i}(1)=\del_{y^{i}},}
which are orthonormal at the respective points $x,z,y$, to identify all of the three tangent spaces canonically with one single Euclidean $\bbR^{n}$. In these coordinates, define an endomorphism $V$ through its matrix representation,
\eq{V=\mrm{diag}\br{1,\fr{1}{2v(1)},\dots,\fr{1}{2v(1)}}.}
 Then
 \eq{\vert \n u_z\vert=\vert V(\n u_y)\vert=\vert V(\n u_{x})\vert\leq \abs{\n u_{x}}=\abs{\n u_{y}},\label{eq:grad norms}}
since 
\eq{v(1)=\fr{1}{2\cos\abs{\dot\G}}\geq \fr 12.}

The second order condition implies that for all $\xi\in \bbR^{n}$ and all $a,b\in \mathrm{GL}_{n}(\bbR)$
\eq{\label{pf:Main-1}0&\leq\xi^{i}\xi^{j}(a_i^p\del_{x^p}+b_i^p\del_{y^p})(a_j^q\del_{x^q}+b_j^q\del_{y^q})Z\\
&=\xi^{i}\xi^{j}\br{a_i^pa_j^q\del_{x^p}\del_{x^q}Z+a_i^pb_j^q\del_{x^p}\del_{y^q}Z+b^{p}_{i}a^{q}_{j}\del_{y^{p}}\del_{x^{q}}Z+b_i^pb_j^q\del_{y^p}\del_{y^q}Z}\\
&=\xi^{i}\xi^{j}a_i^pa_j^q\br{D^2u_z(J_{x^p},J_{x^q})+Du_z(K_{x^px^q})-\fr 12 D^{2}_{pq}u_x}\\
&\hp{=}+\xi^{i}\xi^{j}a_i^pb_j^q(D^2u_z(J_{x^p},J_{y^q})+Du_z(K_{x^py^q}))\\
&\hp{=}+\xi^{i}\xi^{j}b^{p}_i a^{q}_{j}(D^2u_z(J_{y^p},J_{x^q})+Du_z(K_{y^px^q}))\\
&\hp{=}+\xi^{i}\xi^{j}b_i^pb_j^q\br{D^2u_z(J_{y^p},J_{y^q})+Du_z(K_{y^py^q})-\fr 12 D^{2}_{pq}u_y}\\
&=\fr{1}{4}\xi^{i}\xi^{j}(V^{-1}\circ c)_i^p(V^{-1}\circ c)_j^qD^{2}_{pq}u_z-\tfrac{1}{2}(a_i^pa_j^qD^{2}_{pq}u_x+b_i^pb_j^qD^{2}_{pq}u_y)\\
&\hp{=}+\xi^{i}\xi^{j}Du_z(a_i^pa_j^qK_{x^px^q}+a_i^pb_j^q K_{x^py^q}+b^{p}_{i}a^{q}_{j}K_{y^{p}x^{q}}+b_i^pb_j^qK_{y^py^q}),}
where
\eq{
c:=a+b.
}
We test this relation twice with suitable $a$ and $b$ to deduce a concavity relation. First let $\id=a=-b.$
Then we obtain in the sense of bilinear forms on $\bbR^{n}$,
\eq{\label{EC}D^{2}_{ij}u_{x}+D^{2}_{ij}u_{y}\leq 0,}
since in this case $c=0$ and the term involving $K$ also vanishes in view of the $\cK_{3}(0)=0$ case in \Cref{K,,lem:Q=0}. Due to the same reasoning using $\cK_{2}(0)=0$, from setting $a=b=V/2$ we obtain
\eq{-D^{2}_{ij}u_{z}\leq -\fr 12(D^{2}u_{x}+D^{2}u_{y})(V_{i},V_{j}).}
As we are working in Euclidean coordinates, we may view both sides to lie in $\cD_{\bar\G}$. Hence we may apply $f(\abs{\n u_{z}},\cdot)$ after raising an index with respect to the Euclidean scalar product and use the naturality of $f$ to obtain
\eq{b(z,u_{z},\abs{\n u_{z}})&=f(\abs{\n u_{z}},-\n^{2}u_{z})\leq  f\br{\abs{\n u_{z}},-\fr 12 (D^{2}u_{x}+D^{2}u_{y})^{\sharp}},}
where we have used the monotonicity of $f$ in the second variable and \Cref{Courant}.
Using the convexity and the other monotonicity properties, we get
				
\eq{b(z,u_{z},\abs{\n u_{z}})&\leq \fr 12 f(\abs{\n u_{x}},-\n^{2}u_{x})-\fr 12 f(\abs{\n u_{y}},-\n^{2}u_{y})\\
				&=\fr 12 b(x,u_{x},\abs{\n u_{x}})+\fr 12 b(y,u_{y},\abs{\n u_{y}})\\
				&\leq \fr 12 b(x,u_{x},\abs{\n u_{z}})+\fr 12 b(y,u_{y},\abs{\n u_{z}})\\
				&\leq b\br{z,\fr 12 (u_{x}+u_{y}),\abs{\n u_{z}}}.}
Hence, due to the strict monotonicity in the second variable,
\eq{u_{z}\geq \fr 12(u_{x}+u_{y})}
and the proof is complete.

\begin{rem}
The method to introduce different variational direction in the $x$ and $y$ variables was already used in \cite[Thm.~3.13]{Kawohl:/1985}, although in a simplified (quasilinear Euclidean) setting, where it gave a lot more information.
\end{rem}

\section*{Acknowledgments}
This work was made possible through a research scholarship JS received from the DFG and which was carried out at Columbia University in New York. JS would like to thank the DFG, Columbia University and especially Prof.~Simon Brendle for their support.

\bibliographystyle{shamsplain}
\bibliography{Bibliography.bib}
\end{document}